\documentclass[a4paper,draft]{amsart}
\usepackage{amssymb}
\theoremstyle{plain}
  \newtheorem{thm}{Theorem}[section]

  \newtheorem*{obs*}{Observation}
\theoremstyle{definition}

\theoremstyle{remark}
  \newtheorem{rem}[thm]{Remark}
  \newtheorem{com}[thm]{Comment}
  \newtheorem*{ack}{Acknowledgments}
\newcommand{\Z}{\mathbb{Z}}
\newcommand{\C}{\mathbb{C}}
\newcommand{\R}{\mathbb{R}}

\newcommand{\Res}{\operatorname{Res}}
\renewcommand{\Re}{\operatorname{Re}}
\renewcommand{\Im}{\operatorname{Im}}
\numberwithin{equation}{section}
\allowdisplaybreaks
\begin{document}
\title[The colored Jones polynomials of a torus knot]
{Asymptotic behaviors of the colored Jones polynomials of a torus knot}
\author{Hitoshi Murakami}
\address{
Department of Mathematics,
Tokyo Institute of Technology,
Oh-okayama, Meguro, Tokyo 152-8551, Japan
}
\email{starshea@tky3.3web.ne.jp}
\date{\today}
\begin{abstract}
We study the asymptotic behaviors of the colored Jones polynomials
of torus knots.
Contrary to the works by R.~Kashaev, O.~Tirkkonen, Y.~Yokota, and the author,
they do not seem to give the volumes or the Chern--Simons invariants of
the three-manifolds obtained by Dehn surgeries.
On the other hand it is proved that in some cases the limits give the inverse of
the Alexander polynomial.
\end{abstract}
\thanks{This research is partially supported by Grant-in-Aid for Scientific
Research (B) (15340019).}
\maketitle
Let $K$ be a knot in the three-sphere and $J_N(K;t)$ the colored Jones
polynomial corresponding to the $N$-dimensional representation of $sl_2(\C)$
normalized so that $J_N(\text{unknot};t)=1$
\cite{Jones:BULAM385,Kirillov/Reshetikhin:89}.
R.~Kashaev found a series of link invariants parameterized by positive integers
\cite{Kashaev:MODPLA95} and proposed a conjecture that the asymptotic behavior
of his invariants would determine the hyperbolic volume of the knot complement
for any hyperbolic knot \cite{Kashaev:LETMP97}.
It turned out \cite{Murakami/Murakami:ACTAM101} that Kashaev's invariant with
parameter $N$ is equal to
$\left|J_N\left(K;e^{2\pi\sqrt{-1}/N}\right)\right|$.
Kashaev's conjecture was generalized to the volume conjecture which states
that the asymptotic behavior of Kashaev's invariant would determine the
Gromov norm \cite{Gromov:INSHE82} of the knot complement for any knot.
\par
Kashaev and O.~Tirkkonen \cite{Kashaev/Tirkkonen:ZAPNS200} proved the volume
conjecture for torus knots.
More precisely they proved that for the torus knot $T(a,b)$ for coprime integers
$a$ and $b$,
$\displaystyle
\lim_{N\to\infty}
\left.\log{J_N\left(T(a,b),e^{2\pi\sqrt{-1}/N}\right)}\right/N=0$.
Note that since the complement of a torus knot is a Seifert fibered space
\cite{Moser:PACJM1971}, its Gromov norm is zero.
\par
The volume conjecture (or Kashaev's conjecture in this case) is also proved for
the figure-eight knot by T.~Ekholm (see for example
\cite{Murakami:ALDT_VI}).
Moreover in \cite{Murakami/Yokota:2004} Y.~Yokota and the author proved that
for the figure-eight knot $E$ the asymptotic behavior of
$J_N\left(E;e^{2\pi{r}\sqrt{-1}/N}\right)$
determines the volume and the Chern--Simons invariant for the closed
three-manifold obtained by a Dehn surgery along $E$, where $r$ is a complex
parameter near $1$.
More precisely
$\lim_{N\to\infty}\left.\log{J_N\left(E;e^{2\pi{r}\sqrt{-1}/N}\right)}\right/N$
is an analytic function of $r$ that is almost the potential function of
W.~Neumann and D.~Zagier \cite{Neumann/Zagier:TOPOL85} (see also
\cite{Yoshida:INVEM85}).
\par
So one could expect that the colored Jones polynomials of a knot would determine
the volumes and the Chern--Simons invariants for three-manifolds obtained by
Dehn surgeries along any knot.
In this paper, contrary to the results above, we show that the limit of
$\left.\log{J\left(T(a,b);e^{2\pi{r}\sqrt{-1}/N}\right)}\right/N$
behaves strangely; it cannot be extended to a continuous function of $r$.
Note that a three-manifold obtained from a torus knot by Dehn surgery is either
a Seifert fibered space, a lens space, or the connected sum of two lens spaces
\cite{Moser:PACJM1971} and the Gromov norm of such a manifold is $0$.
\begin{ack}
The author thanks K.~Hikami and Y.~Yokota for helpful conversations.
\end{ack}
%%%%%%%%%%%%%%%%%%%%%%%%%%%%%%%%%%%%%%%%%%%%%%%%%%%%%%
\section{Statement of the result}
Let $T(a,b)$ be the $(a,b)$-torus knot for coprime integers $a$ and $b$
($a>1$, $b>1$).
Then its colored Jones polynomial $J_N\bigl(T(a,b);t\bigr)$ is given as follows
\cite{Morton:MATPC95}.
\begin{equation*}
  J_N\bigl(T(a,b);t\bigr)
  =
  \frac{t^{-ab(N^2-1)/4}}{t^{N/2}-t^{-N/2}}
  \sum_{\varepsilon=\pm1}
  \sum_{k=-(N-1)/2}^{(N-1)/2}
  \varepsilon
  t^{abk^2+k(a+\varepsilon{b})+\varepsilon/2}.
\end{equation*}
\par
We will show the following theorem.
%%%%%%%%%%%%%%%%%%%%%%%%%%%%%%%%%%%%%%%%%
\begin{thm}\label{thm}
There exists a neighborhood $U$ of $1\in\C$ such that
if $r\in{U}$ and $r\notin\R$, then
\begin{gather}
  \lim_{N\to\infty}J_N\left(T(a,b);e^{2\pi{r}\sqrt{-1}/N}\right)
  =
  \frac{1}{\Delta\left(T(a,b);e^{2\pi{r}\sqrt{-1}}\right)}
  \quad\text{if $\Im{r}<0$},
  \label{eq:negative}
  \\
  \intertext{and}
  \lim_{N\to\infty}
  \frac{\log{J_N\left(T(a,b);e^{2\pi{r}\sqrt{-1}/N}\right)}}{N}
  =
  \left(1-\dfrac{1}{2abr}-\dfrac{abr}{2}\right)\pi\sqrt{-1}
  \quad\text{if $\Im{r}>0$}.
  \label{eq:positive}
\end{gather}
\end{thm}
\begin{rem}
To be precise, the limit in \eqref{eq:positive} means
\begin{equation*}
  J_N\left(T(a,b);e^{2\pi{r}\sqrt{-1}/N}\right)
  \\
  \underset{N\to\infty}{\sim}
  P(N)e^{N(1-1/(2abr)-abr/2)\pi\sqrt{-1}}
\end{equation*}
for some function of $P(N)$ such that $N^{-c}<|P(N)|<N^c$ for some
positive integer $c$.
Here $f(N)\underset{N\to\infty}{\sim}g(N)$ means that
$\lim_{N\to\infty}f(N)/g(N)=1$.
Note that $\lim_{N\to\infty}\log{P(N)}/N=0$ for such $P(N)$.
Note also that since the real part of the right hand side of \eqref{eq:positive}
is positive, the formula shows that
$\left|J_N\left(T(a,b);e^{2\pi{r}\sqrt{-1}/N}\right)\right|$
grows exponentially if $\Im{r}>0$.
\end{rem}
\begin{rem}
For $r=1$,
Kashaev and Tirkkonen \cite{Kashaev/Tirkkonen:ZAPNS200} proved that
\begin{equation*}
  \left|J_N\left(T(a,b);e^{2\pi{r}\sqrt{-1}/N}\right)\right|
  \underset{N\to\infty}{\sim}
  Ck^{3/2}
\end{equation*}
for some constant $C$.
K.~Hikami and A.~Kirillov study the phase factor \cite{Hikami/Kirillov:PHYLB} and show a relation to the $SU(2)$ Chern--Simons invariant.
\end{rem}
%%%%%%%%%%%%%%%%%%%%%%%%%%%%%%%%%%%%%%%%%%
\section{Proof}
When $t=e^{2\pi{r}\sqrt{-1}/N}$ for $r\in\C\setminus{\Z}$,
$J_N\left(T(a,b);e^{2\pi{r}\sqrt{-1}/N}\right)$ can be given by the following
integral form \cite[Lemma~1]{Kashaev/Tirkkonen:ZAPNS200}.
\begin{equation*}
  J_N\left(T(a,b),e^{2\pi{r}\sqrt{-1}/N}\right)
  =
  \Phi_{a,b,r}(N)
  \int_{C}
  e^{Nf_{a,b,r}(z)}
  \tau_{a,b}(z)
  dz.
\end{equation*}
Here
\begin{align*}
  \Phi_{a,b,r}(N)
  &:=g_{a,b,r}\sqrt{N}
     e^{-\left(ab\left(N^2-1\right)+a/b+b/a\right)\pi{r}\sqrt{-1}/(2N)}
  \\
  &\text{
    with
      $g_{a,b,r}
       :=
       \frac{\sqrt{ab}}
            {2\pi\sqrt{2r}
            e^{\pi\sqrt{-1}/4}\sinh\left(\pi{r}\sqrt{-1}\right)}$},
  \\[3mm]
  f_{a,b,r}(z)
  &:=
  ab\left(z-\frac{z^2}{2\pi{r}\sqrt{-1}}\right),
  \\[3mm]
  \tau_{a,b}(z)
  &:=
  \frac{2\sinh(az)\sinh(bz)}{\sinh(abz)},
\end{align*}
and $C$ is the line
\begin{equation}\label{eq:C}
  C:=\{se^{\varphi\sqrt{-1}}\mid{s}\in\R\}
\end{equation}
with $\Re\left(r\sqrt{-1}e^{-2\varphi\sqrt{-1}}\right)>0$.
(We will choose $\varphi$ later.)
\par
Let $\mathcal{P}$ be the set of the poles of $\tau_{a,b}(z)$, that is , we put
\begin{equation*}
  \mathcal{P}
  :=
  \left\{\frac{k\pi\sqrt{-1}}{ab} \big| k\in\Z,a\nmid k,b\nmid k\right\}.
\end{equation*}
We will show that the result follows if $r\notin\R$, $\Re{r}>0$,
and $|r|>\frac{1}{ab}$.
\par
Put $\theta:=\arg{r}$ ($-\pi/2<\theta<\pi/2$) so that
$r=|r|e^{\theta\sqrt{-1}}$.
Put also $\varphi:=\theta/2+\pi/4$ in \eqref{eq:C}.
Note that $\Re\left(r\sqrt{-1}e^{-2\varphi\sqrt{-1}}\right)=\Re|r|>0$
and that $0<\varphi<\pi/2$.
\par
Let ${C'}$ be the line parallel to $C$ that passes through the point
$\pi{r}\sqrt{-1}$.
(We change $C'$ slightly near $\mathcal{P}$ to avoid poles of $\tau_{a,b}(z)$
if necessary.)
For a positive number $R$, let $D_{+}$ ($D_{-}$, respectively) be the segment of
the line $\Re{z}=R$ ($\Re{z}=-R$, respectively) between $C$ and $C'$ oriented
upward (downward, respectively).
We see that $C'$ and $D_{\pm}$ are parameterized as follows.
\begin{equation*}
\begin{split}
  C'
  &=
  \left\{
    se^{(\theta/2+\pi/4)\sqrt{-1}}+\pi{r}\sqrt{-1}
    \Bigm|
    s\in\R
  \right\},
  \\
  D_{\pm}
  &=
  \left\{
    \vphantom{\bigr)}
    \pm{R}+s\sqrt{-1}
    \bigm|
    \pm{R}\tan\left(\theta/2+\pi/4\right)
    \le{s}
    \le
    \pm{R}\tan\left(\theta/2+\pi/4\right)
    +
    \pi|r|h(\theta)
  \right\},
\end{split}
\end{equation*}
where $h(\theta):=\cos\theta+\sin\theta\tan(\theta/2+\pi/4)$
Note that $C'$ crosses the imaginary axis at $\pi|r|h(\theta)\sqrt{-1}$.
\par
Then by the residue theorem
\begin{equation*}
\begin{split}
  &
  \int_{\overline{C}}
  e^{Nf_{a,b,r}(z)}\tau_{a,b}(z)\,dz
  \\
  &\quad
  =
  \int_{\overline{C'}}
  e^{Nf_{a,b,r(z)}}\tau_{a,b}(z)\,dz
  +
  \sum_{k}
  \Res
  \left(
    e^{Nf_{a,b,r}(z)}\tau_{a,b}(z);z=k\pi\sqrt{-1}/(ab)
  \right)
  \\
  &\quad\quad
  -\int_{D_+}
  e^{Nf_{a,b,r}(z)}\tau_{a,b}(z)\,dz
  -\int_{D_-}
  e^{Nf_{a,b,r}(z)}\tau_{a,b}(z)\,dz,
\end{split}
\end{equation*}
where $\Res(F(z);z=\zeta)$ is the residue of $F(z)$ around $z=\zeta$,
$\overline{C}$ ($\overline{C'}$, respectively)
is the segment in $C$ (${C'}$, respectively) bounded by
$D_{\pm}$ oriented from left to right,
and $k$ runs over integers that are not multiples of $a$ or $b$ such that
$k\pi\sqrt{-1}/ab$ is between $C$ and $C'$, that is, $0<k<ab|r|h(\theta)$.
\par
Since it is not hard to see that on $D_{\pm}$
\begin{equation*}
  \left|e^{Nf_{a,b,r}(z)}\right|
  <
  e^{abN\left(c_2R^2+c_1R+c_0\right)}
\end{equation*}
for some $c_i$ with $c_2<0$ if $R$ is sufficiently large, and that
\begin{equation*}
\begin{split}
  \left|\tau_{a,b}(z)\right|
  &=
  \frac{\left|\left(e^{az}-e^{-az}\right)
              \left(e^{bz}-e^{-bz}\right)\right|}
       {\left|e^{abz}-e^{-abz}\right|}
  \\
  &\le
  \frac{\left|e^{(a+b)z}\right|
        +
        \left|e^{(a-b)z}\right|
        +
        \left|e^{(b-a)z}\right|
        +
        \left|e^{-(a+b)z}\right|}
       {\bigl|\left|e^{abz}\right|-\left|e^{-abz}\right|\bigr|}
  \\
  &=
  \frac{e^{\pm(a+b)R}
        +
        e^{\pm(a-b)R}
        +
        e^{\pm(b-a)R}
        +
        e^{\mp(a+b)R}}
       {\left|e^{\pm{abR}}-e^{\mp{abR}}\right|},
\end{split}
\end{equation*}
the integrals on $D_{\pm}$ vanish when $R\to\infty$.
Therefore we have
\begin{multline*}
  \int_{C}
  e^{Nf_{a,b,r}(z)}\tau_{a,b}(z)\,dz
  \\
  =
  \int_{C'}
  e^{Nf_{a,b,r}(z)}\tau_{a,b}(z)\,dz
  +
  \sum_{k}
  \Res
  \left(
    e^{Nf_{a,b,r}(z)}\tau_{a,b}(z);z=k\pi\sqrt{-1}/(ab)
  \right).
\end{multline*}
\par
We will apply the saddle point method (or the method of steepest descent)
to know the asymptotic behavior of the integral on $C'$.
See for example \cite[Theorem 7.2.9]{Marsden/Hoffman:Complex_Analysis}.
On $C'$ we have
\begin{equation*}
  f_{a,b,r}\left(se^{(\theta/2+\pi/4)\sqrt{-1}}+\pi{r}\sqrt{-1}\right)
  =
  -\frac{s^2}{2\pi|r|}+\frac{\pi{r}\sqrt{-1}}{2}.
\end{equation*}
So $\Im{f_{a,b,r}(z)}$ is constant on $C'$ and $\Re{f_{a,b,r}(z)}$ takes its
unique maximum at $z=\pi{r}\sqrt{-1}$ (when $s=0$) on $C'$.
(Note that this is also true if we change $C'$ slightly near a pole of
$\tau_{a,b}$.)
Therefore from the saddle point method we have
\begin{equation*}
\begin{split}
  &
  \int_{C'}
  e^{Nf_{a,b,r}(z)}\tau_{a,b}(z)\,dz
  \\
  &\quad\underset{N\to\infty}{\sim}
  \frac{\sqrt{2\pi}\tau_{a,b}(\pi{r}\sqrt{-1})}
       {\sqrt{N}\sqrt{-d^2f_{a,b,r}(z)/dz^2\big|_{z=\pi{r}\sqrt{-1}}}}
  e^{Nf_{a,b,r}(\pi{r}\sqrt{-1})}
  \\
  &=
  \pi\sqrt{\frac{2r}{abN}}e^{\pi\sqrt{-1}/4}
  \tau_{a,b}(\pi{r}\sqrt{-1})
  e^{Nab\pi{r}\sqrt{-1}/2}.
\end{split}
\end{equation*}
We also have
\begin{equation*}
\begin{split}
  &
  \Res
  \left(
    e^{Nf_{a,b,r}(z)}\tau_{a,b}(z);z=k\pi\sqrt{-1}/(ab)
  \right)
  \\
  &\quad=
  e^{Nk\pi\sqrt{-1}\left(1-\frac{k}{2abr}\right)}
  \Res\left(\tau_{a,b}(z);z=k\pi\sqrt{-1}/(ab)\right)
  \\
  &\quad=
  (-1)^k
  \frac{2\sinh(k\pi\sqrt{-1}/a)\sinh(k\pi\sqrt{-1}/b)}{ab}
  e^{Nk\pi\sqrt{-1}\left(1-\frac{k}{2abr}\right)}.
\end{split}
\end{equation*}
\par
Therefore we have
\begin{equation}\label{eq:asymptotic}
\begin{split}
  &\int_{C}
  e^{Nf_{a,b,r}(z)}\tau_{a,b}(z)\,dz
  \\
  &\quad\underset{N\to\infty}{\sim}
  \pi\sqrt{\frac{2r}{abN}}e^{\pi\sqrt{-1}/4}
  \tau_{a,b}(\pi{r}\sqrt{-1})
  e^{Nab\pi{r}\sqrt{-1}/2}
  \\
  &\quad\phantom{\underset{N\to\infty}{\sim}}
  +
  \sum_{k}
  (-1)^k
  \frac{2\sinh(k\pi\sqrt{-1}/a)\sinh(k\pi\sqrt{-1}/b)}{ab}
  e^{Nk\pi\sqrt{-1}\left(1-\frac{k}{2abr}\right)}.
\end{split}
\end{equation}
\par
Now we want to know the biggest real part of the exponents in
\eqref{eq:asymptotic}.
We have
\begin{equation*}
\begin{split}
  &
  \max
  \left\{
    \Re\left(\frac{ab\pi{r}\sqrt{-1}}{2}\right),
    \max_{k}{\Re\left(k\pi\sqrt{-1}\left(1-\frac{k}{2abr}\right)\right)}
  \right\}
  \\
  &\quad=
  \max
  \left\{
    -\frac{ab\pi\Im{r}}{2},
    \max_{k}{\frac{-k^2\pi\Im{r}}{2ab|r|^2}}
  \right\}
  \\
  &\quad=
  \begin{cases}
    \displaystyle
    \frac{\pi\Im{r}}{2ab}\max_{k}
    \left\{
      -a^2b^2,-\frac{k^2}{|r|^2}
    \right\}
    &\qquad\text{if $\Im{r}>0$}
    \\[5mm]
    \displaystyle
    -\frac{\pi\Im{r}}{2ab}\max_{k}
    \left\{
      a^2b^2,\frac{k^2}{|r|^2}
    \right\}
    &\qquad\text{if $\Im{r}<0$}
  \end{cases}.
\end{split}
\end{equation*}
Since $1<h(\theta)$ when $0<\theta<\pi/2$, $|r|>1/(ab)$, and $k$
runs over positive integers less than $ab|r|h(\theta)>1$, we have
\begin{equation*}
  \max_{k}\left\{-a^2b^2,-\frac{k^2}{|r|^2}\right\}
  =
  -\frac{1}{|r|^2}
\end{equation*}
if $\Im{r}>0$.
Since $0<h(\theta)<1$ when $-\pi/2<\theta<0$, $0<k<ab|r|h(\theta)<ab|r|$.
So we have
\begin{equation*}
  \max_{k}\left\{a^2b^2,\frac{k^2}{|r|^2}\right\}
  =
  a^2b^2
\end{equation*}
if $\Im{r}<0$.
\par
Therefore if $\Im{r}<0$
\begin{equation*}
  \int_{C}
  e^{Nf_{a,b,r(z)}}\tau_{a,b}(z)dz
  \underset{N\to\infty}{\sim}
  \pi\sqrt{\frac{2r}{abN}}e^{\pi\sqrt{-1}/4}
  \tau_{a,b}(\pi{r}\sqrt{-1})
  e^{Nab\pi{r}\sqrt{-1}/2}
\end{equation*}
and if $\Im{r}>0$
\begin{equation*}
  \int_{C}
  e^{Nf_{a,b,r}(z)}\tau_{a,b}(z)\,dz
  \\
  \underset{N\to\infty}{\sim}
  -
  \frac{2\sinh(\pi\sqrt{-1}/a)\sinh(\pi\sqrt{-1}/b)}{ab}
  e^{N\pi\sqrt{-1}\left(1-\frac{1}{2abr}\right)}.
\end{equation*}
Since
\begin{equation*}
  \Phi_{a,b,r}(N)
  \underset{N\to\infty}{\sim}
  g_{a,b,r}\sqrt{N}e^{-Nabr\pi\sqrt{-1}/2},
\end{equation*}
we have
\begin{multline*}
  J_N\bigr(T(a,b);e^{2\pi{r}\sqrt{-1}}\bigr)
  \\
  \underset{N\to\infty}{\sim}
  -
  \frac{2\sinh(\pi\sqrt{-1}/a)\sinh(\pi\sqrt{-1}/b)g_{a,b,r}}{ab}
  \sqrt{N}
  e^{N\pi\sqrt{-1}\bigl(1-abr/2-1/(2abr)\bigr)}
\end{multline*}
if $\Im{r}>0$ and
\begin{equation*}
  J_N\bigr(T(a,b);e^{2\pi{r}\sqrt{-1}}\bigr)
  \underset{N\to\infty}{\sim}
  \frac{\sinh(ar\pi\sqrt{-1})\sinh(br\pi\sqrt{-1})}
       {\sinh(abr\pi\sqrt{-1})\sinh(r\pi\sqrt{-1})}
\end{equation*}
if $\Im{r}<0$.
\par
Since the Alexander polynomial $\Delta\bigl(T(a,b);t\bigr)$ of $T(a,b)$
is
$\frac{(t^{ab/2}-t^{-ab/2})(t^{1/2}-t^{-1/2})}
      {(t^{a/2}-t^{-a/2})(t^{b/2}-t^{-b/2})}$
(see for example \cite[Page 119]{Lickorish:1997}), we have
\begin{equation*}
  \lim_{N\to\infty}J_N\left(T(a,b);e^{2\pi{r}\sqrt{-1}/N}\right)
  =
  \frac{1}{\Delta\left(T(a,b);e^{2\pi{r}\sqrt{-1}}\right)}
\end{equation*}
if $\Im{r}<0$.
If $\Im{r}>0$ we have
\begin{equation*}
  \lim_{N\to\infty}
  \frac{\log{J_N\left(T(a,b);e^{2\pi{r}\sqrt{-1}/N}\right)}}{N}
  =
  \left(1-\dfrac{1}{2abr}-\dfrac{abr}{2}\right)\pi\sqrt{-1}.
\end{equation*}
This completes the proof of Theorem~\ref{thm}.
%%%%%%%%%%%%%%%%%%%%%%%%%%%%%%%%%%%%%%%%%%%%%%%%%%%%%%%%%%%%%%%%%%%%%
\section{Comments}
Here are a couple of probably inaccurate comments about the relation of the
limit of the colored Jones polynomials of a torus knot and its Alexander
polynomial.
\begin{com}
The Melvin--Morton--Rozansky conjecture proved by D.~Bar-Natan and
S.~Garoufalidis \cite{BarNatan/Garoufalidis:INVEM96}, which was proposed by
P.~Melvin and H.~Morton \cite{Melvin/Morton:COMMP95} and `proved' non-rigorously
by L.~Rozansky \cite{Rozansky:COMMP96}, tells that if we expand the colored
Jones polynomial of a knot as a power series
\begin{equation*}
  J_N\left(K;e^{h}\right)=\sum_{d,l\ge0}b_{dl}(K)h^dN^l,
\end{equation*}
then
\begin{enumerate}
  \item[(i)]
  $b_{dl}(K)=0$ if $l>d$,
  \item[(ii)]
  the diagonal coefficients give the inverse of the Alexander polynomial,
  that is,
  \begin{equation*}
    \sum_{d\ge0}b_{dd}(K)(hN)^d=\frac{1}{\Delta(K;e^{hN})},
  \end{equation*}
  where $\Delta(K;t)$ is the Alexander polynomial of $K$.
\end{enumerate}
If we replace $h$ with $2\pi{r}\sqrt{-1}/N$, we have from (i)
\begin{equation}\label{eq:MMR1}
  J_N\left(K;e^{2\pi{r}\sqrt{-1}/N}\right)
  =
  \sum_{d\ge{l}\ge0}b_{dl}(K)(2\pi{r}\sqrt{-1})^dN^{l-d}
\end{equation}
and from (ii)
\begin{equation}\label{eq:MMR2}
  \sum_{d\ge0}b_{dd}(K)(2\pi{r}\sqrt{-1})^d
  =
  \frac{1}{\Delta(K;e^{2\pi{r}\sqrt{-1}})}.
\end{equation}
\par
Theorem~\ref{thm} would say that if $\Im{r}<0$ and $K$ is a torus knot,
then the right hand side of \eqref{eq:MMR1} converges to the left hand side
of \eqref{eq:MMR2} and gives the inverse of the Alexander polynomial.
\par
Note that this does not hold for the figure-eight knot
\cite{Murakami/Yokota:2004}.
\end{com}
\begin{com}
In \cite[Theorem 4]{Hikami:2004} Hikami obtains a recursive formula for
the colored Jones polynomials of a torus knot:
\begin{multline}\label{eq:recursive}
  J_N\bigl(T(a,b);t\bigr)
  \\
  =
  \frac{t^{(a-1)(b-1)(1-N)/2}}{1-t^{-N}}
  \left(
    1-t^{a(1-N)-1}-t^{b(1-N)-1}+t^{(a+b)(1-N)}
  \right)
  \\
  +
  \frac{1-t^{2-N}}{1-t^{-N}}
  t^{ab(1-N)-1}
  J_{N-2}\bigl(T(a,b),t\bigr).
\end{multline}
\par
The limit \eqref{eq:negative} could be obtained from the following `fake'
calculation.
If we replace $t$ with $e^{2\pi{r}\sqrt{-1}/N}$
in \eqref{eq:recursive}, we could approximate
\begin{equation*}
  t^{(a-1)(b-1)(1-N)/2}
  =
  e^{r(a-1)(b-1)\pi\sqrt{-1}/N-r(a-1)(b-1)\pi\sqrt{-1}}
\end{equation*}
as $e^{-r(a-1)(b-1)\pi\sqrt{-1}}$ for large $N$.
Approximating the other terms similarly, we would have
\begin{equation*}
\begin{split}
  J_{\infty}
  &=
  \frac{e^{-r(a-1)(b-1)\pi\sqrt{-1}}}{1-e^{-2r\pi\sqrt{-1}}}
  \left(
    1-e^{-2ar\pi\sqrt{-1}}-e^{-2br\pi\sqrt{-1}}+e^{-2(a+b)r\pi\sqrt{-1}}
  \right)
  \\
  &\quad+
  e^{-2abr\pi\sqrt{-1}}J_{\infty},
\end{split}
\end{equation*}
where $J_{\infty}$ denotes the `limit' of
$J_N\left(T(a,b);e^{2\pi{r}\sqrt{-1}/N}\right)$.
Then \eqref{eq:negative} follows easily.
\par
Note that there exists a recursive formula for the colored Jones polynomials of
any knot \cite{Garoufalidis/Le:2003}.
See also \cite{Gelca:PROAM2002,Gelca/Sain:JKNOT2003} for a recursive formula
of the $(2,2n+1)$ torus knot.
\end{com}
%%%%%%%%%%%%%%%%%%%%%%%%%%%%%%%%%%%%%%%%%%%%%
\bibliography{mrabbrev,hitoshi}

\providecommand{\bysame}{\leavevmode\hbox to3em{\hrulefill}\thinspace}
\begin{thebibliography}{10}

\bibitem{BarNatan/Garoufalidis:INVEM96}
D.~Bar-Natan and S.~Garoufalidis, \emph{On the {M}elvin-{M}orton-{R}ozansky
  conjecture}, Invent. Math. \textbf{125} (1996), no.~1, 103--133.

\bibitem{Garoufalidis/Le:2003}
S.~Garoufalidis and T.~T.~Q. Le, \emph{{The colored Jones function is
  q-holonomic}}, \mbox{arXiv:math.GT/0309214}.

\bibitem{Gelca:PROAM2002}
R.~Gelca, \emph{On the relation between the {$A$}-polynomial and the {J}ones
  polynomial}, Proc. Amer. Math. Soc. \textbf{130} (2002), no.~4, 1235--1241
  (electronic).

\bibitem{Gelca/Sain:JKNOT2003}
R.~Gelca and J.~Sain, \emph{The noncommutative {A}-ideal of a
  {$(2,2p+1)$}-torus knot determines its {J}ones polynomial}, J. Knot Theory
  Ramifications \textbf{12} (2003), no.~2, 187--201.

\bibitem{Gromov:INSHE82}
M.~Gromov, \emph{Volume and bounded cohomology}, Inst. Hautes \'Etudes Sci.
  Publ. Math. (1982), no.~56, 5--99 (1983).

\bibitem{Hikami:2004}
K.~Hikami, \emph{{Difference equation of the colored Jones polynomial for torus
  knot}}, \mbox{arXiv:math.GT/0403224}.

\bibitem{Hikami/Kirillov:PHYLB}
K.~Hikami and A.~N. Kirillov, \emph{Torus knot and minimal model}, Phys. Lett.
  B \textbf{575} (2003), 343--348.

\bibitem{Jones:BULAM385}
V.~F.~R. Jones, \emph{A polynomial invariant for knots via von {N}eumann
  algebras}, Bull. Amer. Math. Soc. (N.S.) \textbf{12} (1985), no.~1, 103--111.

\bibitem{Kashaev:MODPLA95}
R.~M. Kashaev, \emph{A link invariant from quantum dilogarithm}, Modern Phys.
  Lett. A \textbf{10} (1995), no.~19, 1409--1418.

\bibitem{Kashaev:LETMP97}
\bysame, \emph{The hyperbolic volume of knots from the quantum dilogarithm},
  Lett. Math. Phys. \textbf{39} (1997), no.~3, 269--275.

\bibitem{Kashaev/Tirkkonen:ZAPNS200}
R.~M. Kashaev and O.~Tirkkonen, \emph{A proof of the volume conjecture on torus
  knots}, Zap. Nauchn. Sem. S.-Peterburg. Otdel. Mat. Inst. Steklov. (POMI)
  \textbf{269} (2000), no.~Vopr. Kvant. Teor. Polya i Stat. Fiz. 16, 262--268,
  370.

\bibitem{Kirillov/Reshetikhin:89}
A.~N. Kirillov and N.~Yu. Reshetikhin, \emph{Representations of the algebra
  ${U_q(sl(2))}$, $q$-orthogonal polynomials and invariants of links}, Infinite
  Dimensional Lie Algebras and Groups (V.G. Kac, ed.), Advanced {S}eries in
  {M}athematical {P}hysics, vol.~7, World Scientifics, Singapore, 1989.

\bibitem{Lickorish:1997}
W.~B.~R. Lickorish, \emph{An introduction to knot theory}, Graduate Texts in
  Mathematics, vol. 175, Springer-Verlag, New York, 1997.

\bibitem{Marsden/Hoffman:Complex_Analysis}
J.~E. Marsden and M.~J. Hoffman, \emph{Basic complex analysis}, W. H. Freeman
  and Company, New York, 1987.

\bibitem{Melvin/Morton:COMMP95}
P.~M. Melvin and H.~R. Morton, \emph{The coloured {J}ones function}, Comm.
  Math. Phys. \textbf{169} (1995), no.~3, 501--520.

\bibitem{Morton:MATPC95}
H.~R. Morton, \emph{The coloured {J}ones function and {A}lexander polynomial
  for torus knots}, Math. Proc. Cambridge Philos. Soc. \textbf{117} (1995),
  no.~1, 129--135.

\bibitem{Moser:PACJM1971}
L.~Moser, \emph{Elementary surgery along a torus knot}, Pacific J. Math.
  \textbf{38} (1971), 737--745.

\bibitem{Murakami:ALDT_VI}
H.~Murakami, \emph{The asymptotic behavior of the colored {J}ones function of a
  knot and its volume}, Proceedings of `Art of Low Dimensional Topology VI'
  (T.~Kohno, ed.), January 2000, pp.~87--96, \mbox{arXiv:math.GT/0004036}.

\bibitem{Murakami/Murakami:ACTAM101}
H.~Murakami and J.~Murakami, \emph{The colored {J}ones polynomials and the
  simplicial volume of a knot}, Acta Math. \textbf{186} (2001), no.~1, 85--104.

\bibitem{Murakami/Yokota:2004}
H.~Murakami and Y.~Yokota, \emph{{The colored Jones polynomials of the
  figure-eight knot and its Dehn surgery spaces}},
  \mbox{arXiv:math.GT/0401084}.

\bibitem{Neumann/Zagier:TOPOL85}
W.~D. Neumann and D.~Zagier, \emph{Volumes of hyperbolic three-manifolds},
  Topology \textbf{24} (1985), no.~3, 307--332.

\bibitem{Rozansky:COMMP96}
L.~Rozansky, \emph{A contribution of the trivial connection to the {J}ones
  polynomial and {W}itten's invariant of $3$d manifolds. {I}, {I}{I}}, Comm.
  Math. Phys. \textbf{175} (1996), no.~2, 275--296, 297--318.

\bibitem{Yoshida:INVEM85}
T.~Yoshida, \emph{The $\eta$-invariant of hyperbolic $3$-manifolds}, Invent.
  Math. \textbf{81} (1985), no.~3, 473--514.

\end{thebibliography}
\bibliographystyle{hamsplain}
\end{document}